\documentclass[oneside,12pt,a4paper]{article}

\usepackage{mathptmx}
\usepackage[english]{babel}
\usepackage[fleqn]{amsmath}
\usepackage{amsfonts}
\usepackage{dsfont}
\usepackage{amsxtra,amscd}
\usepackage{latexsym}
\usepackage{amssymb}
\usepackage{enumerate}
\usepackage{eucal}
\usepackage{graphicx}
\usepackage{psfrag}
\usepackage{epsfig}
\usepackage[colorlinks = true,
            linkcolor = red,
            urlcolor  = blue,
            citecolor = blue,
            anchorcolor = blue]{hyperref}
\usepackage{cite}

\usepackage[misc]{ifsym} 

\usepackage[left=2.5cm,right=2.5cm,top=2cm,bottom=2cm,includeheadfoot]{geometry}

\newtheorem{theo}{Theorem}
\newtheorem{Def}{Definition}
\newtheorem{Lem}{Lemma}

\newtheorem{cor}{Corollary}
\newcommand{\rem}{{\bf Remark}}
\newcommand{\bsp}{{\bf Example\,}}
\newcommand{\la}{\lambda}
\newcommand{\prf}{\textit{Proof.\ }}
\newcommand{\hs}{\,}
\newcommand{\ou}{$\mathbf{1}$}
\newcommand{\be}{\begin{equation}}
\newcommand{\ee}{\end{equation}}
\newcommand{\Sup}{\vee}
\newcommand{\Inf}{\wedge}
\newcommand{\nl}{\left\|}          
\newcommand{\nr}{\right\|}           
\newcommand{\bl}{\left|}
\newcommand{\br}{\right|}
\newcommand{\Lra}{\Longrightarrow}
\newcommand{\Llra}{\Longleftrightarrow}
\newcommand{\Ra}{\Rightarrow}
\newcommand{\La}{\Leftarrow}
\newcommand{\Orth}{\rm Orth}
\newcommand{\bR}{{\mathbb R}}
\newcommand{\bN}{\mathbb N}
\newcommand{\cA}{\mathcal A}
\newcommand{\cP}{\mathcal P}
\newcommand{\fP}{\mathfrak P}

\newcommand\qed{\hfill\quad $\Box$}

\psfrag{phi}[]{\small{$\varphi$}}
\psfrag{pi}[]{\small{$\pi$}}
\psfrag{tilf}[]{\small{$\tilde{\varphi}$}}
\psfrag{f1}[]{\small{$f_1$}}
\psfrag{f2}[]{\small{$f_2$}}
\psfrag{f3}[]{\small{$f_3$}}

\begin{document}
\title{On Finite Elements in $f$-Algebras and in Product Algebras}
\author{Helena Malinowski\footnote{H. Malinowski, Technische Universit\"at Dresden, Germany \hspace*{2mm}\Letter\hspace*{2mm}lenamalinowski@gmx.de} \and Martin~R.~Weber\footnote{M.~R.~Weber, Technische Universit\"at Dresden, Germany \hspace*{2mm}\Letter\hspace*{2mm}martin.weber@tu-dresden.de}}

\date{\vspace{-5ex}}
\maketitle

\begin{center}
26 September 2012
\end{center}
\begin{abstract}
Finite elements, which are well-known and studied in the framework of vector lattices, 
are investigated in $\ell$-algebras, preferably in $f$-algebras, and in product algebras.  
The additional structure of an associative multiplication leads to new questions and some new properties concerning the collections of finite, totally finite and self-majorizing elements. 
In some cases the order ideal of finite elements is a ring ideal as well. 
It turns out that a product of elements in an $f$-algebra is a finite element if at least one factor is finite. If the multiplicative unit exists, the latter plays an important role in the investigation of finite elements. For the product of special $f$-algebras an element is finite in the algebra if and only if its power is finite in the product algebra.
\end{abstract}
\section{Introduction}\label{Sec1}
Finite  elements in Archimedean vector lattices were introduced 1971-1974 in the papers \cite{MakWeb74,MakWeb77a,Web73} as an abstract analogue of continuous functions (on a locally compact space) with compact support. 
Finite and totally finite elements play a very important role in the representation theory of Archimedean 
vector lattices by means of real valued (i.e. everywhere finite) continuous functions on a locally compact Hausdorff space, where they are required to be represented as functions with compact support. 
The classes of these elements in general vector lattices and Banach lattices are thoroughly studied in 
a number of papers, see \cite{ChenWeb03-1,ChenWeb03-2,MakWeb74,MakWeb77a,MakWeb77b,MakWeb78,Web95b}. Finite elements in vector lattices of operators are dealt with in \cite{ChenWeb03-3,HHW}. 
A condensed short overview concerning finite and totally finite elements the reader can find in \cite{Web06}.
Self-majorizing elements in Archimedean vector lattices (also known as semi-order units) have been studied  systematically in the recent paper \cite{TeiWeb11}. 
\par\smallskip
In this paper we investigate finite elements in Archimedean $\ell$-, $d$- and $f$-algebras and in pro\-duct algebras. 
It is well known that the  vector lattice of all orthomorphisms on an Archimedean vector lattice is an Archimedean $f$-algebra with weak order unit, see e.g. \cite[Theorem 8.24]{AliBur85}.  We use this fact in Sections \ref{Sec3} and \ref{Sec5}.
\par\medskip 
The material of the paper is organized as follows: in Section \ref{Sec2} we provide the notions of the theory of vector lattices and $\ell$-algebras which are necessary in order to present our results. In particular, 
we define the finite, totally finite and self-majorizing elements in an Archimedean vector lattice. Further we list some properties of $f$-algebras which are relevant for our purpose.
In Section \ref{Sec3} we study finite elements in $f$- and $d$-algebras with multiplicative unit, in Section \ref{Sec4} we investigate them in $f$-algebras without multiplicative unit. For this purpose the weak factorization pro\-per\-ty is introduced and its relations to other well-known properties in $f$-algebras are demonstrated by examples. In Section \ref{Sec5} we consider finite elements in products of uniformly complete $f$-algebras.
For details concerning vector lattices and $\ell$-algebras we refer to the monographs \cite{AbrAli02,AliBur85,Mey91,Zaa83,Zaa97} as well as to the papers \cite{BerHui90,BeuHui84,Bou00,Bou03,Hui91,deP81}. The recent development in the theory of $\ell$-algebras is reflected in the survey paper \cite{BBT03}.
\section{Preliminaries}\label{Sec2}
Recall some definitions and notations known in the theory of vector lattices and lattice ordered algebras. 
A vector lattice will be denoted by $V$, a Banach lattice by $E$ and, a lattice ordered algebra by $\cA$. 
By $V_+, E_+$ and $\cA_+$ will be denoted their cones of positive elements, respectively.
We consider only \textbf{Archimedean} vector lattices $V$ and $\cA$. This assumption, in particular, ensures 
the uniqueness of uniform limits (see \cite{LuxZaa71}).  
For details see \cite{AbrAli02,AliBur85,BerHui90,Mey91,Zaa83}.
\begin{itemize}
\item If $A$ is a non-empty subset of $V$ then the smallest order ideal that contains $A$ is denoted by $I_{A}$ and 
      is called the {\em order ideal generated by $A$}. This order ideal is (see \cite{AliBur85})
      \[ I_A=\left\{x\in V\colon \exists a_1,\ldots,a_n\in A \mbox{ and } \la_1,\ldots,\la_n\in \bR^+ 
       \mbox{ such that } \bl x\br\leq\sum_{i=1}^n\la_i\bl a_i\br \right\}. \]
\item For a non-empty subset $A\subseteq V$ by $A^\perp$ we denote the set
      $\left\{x\in V \colon \forall a\in A\hspace*{1mm} x\perp a\right\}$.
\item The set $A^{\perp\perp}$ is known as \textit{the band generated by} $A$, i.e. the smallest band that contains $A$.
      If $A$ consists of one element $x$, then $\left\{x\right\}^{\perp\perp}$ is called the \textit{principal band} generated by $x$.
\item An element $u\in V_+$ is an (strong) {\it order unit}, if for each $x\in V$ there 
      is a $\la\in \bR_{\geq 0}$ with $-\la u\leq x\leq \la u$ (or equivalently, 
      $|x|\leq \la u $).
\item An element $u\in V_+$ is a {\it weak order unit}, if $\{u\}^{\perp\perp}=V$, i.e. $x\in V$ and 
      $x\perp u$ imply $x=0$. 
\end{itemize}
\medskip
Further on an \textit{algebra} is understood to be a set $\cA$ equipped with several operations: 
beside the addition ($+$) and the usual scalar multiplication, which turn $\cA$ into a vector space, there is also defined an associative multiplication ($\cdot$) satisfying the distributive laws.
\begin{itemize}
\item A vector lattice $\cA$ is called a \textit{lattice ordered algebra}, a \textit{Riesz-algebra} 
      or also an $\ell$-\textit{algebra}, if $\cA$ is equipped with an associative 
      multiplication\footnote{It is convenient to write $ab$ instead of $a\cdot b$ for the product of $a$ 
      and $b$.} such that $\cA$ becomes an algebra, where 
      \be\label{ell} (\ell) \qquad a,b\geq 0 \hs\Lra \hs ab\geq 0 \;\mbox{ holds for all } \, 0\leq a,b\in \cA. 
      \nonumber\ee 
The basic notions and properties of $\ell$-algebras can be found in \cite[Chapter 20]{Zaa83}.
      Equivalent to ($\ell$) are the conditions: \\
      \hspace*{1cm}  
                   (${\ell_1}$) \; if $a,b,c \in \cA$ satisfy $a\leq b$ and $c\geq 0$ then 
                   $ac\leq bc$,   \\
      \hspace*{1cm}
                   (${\ell_2}$) \; $|ab|\leq |a||b|\,$ for all $a,b\in \cA$, \\ 
      see \cite[Sect.~1]{BerHui90}.  
\item An $\ell$-algebra is called a $d$-\textit{algebra} (see \cite{Kud62}), if it satisfies the condition 
      \be\label{d} ($d$) \qquad a\Inf b=0  \hs\Lra \hs (ac)\Inf (bc) = (ca)\Inf(cb)=0 \mbox{ for all } 
      c\geq 0. \nonumber \ee 
      Equivalent to (${\rm d}$) are the conditions: \\ 
      \hspace*{1cm}  
                   (${\rm d_1}$) \; $|ab|=|a||b|$ for all $a,b\in \cA$, and also \\
      \hspace*{1cm}
                   (${\rm d_2}$) \; $c(a\Inf b) = ca \Inf cb \, \mbox{ and }\, (a\Inf b)c = ac \Inf bc 
       \mbox{ for all }a,b\in \cA, c\in \cA_+$, \\ 
       \hspace*{1cm}
                   (${\rm d_3}$) \; $c(a\Sup b) = ca \Sup cb \, \mbox{ and }\, (a\Sup b)c = ac \Sup bc 
       \mbox{ for all }a,b\in \cA, c\in \cA_+$, \\ 
       see \cite[Proposition 1.2]{BerHui90}. 
\item An $\ell$-algebra is called an $f$-\textit{algebra}, if it satisfies the condition 
      \be\label{f} ($f$) \qquad a\Inf b=0  \mbox{ for all } c\geq 0 \hs\Lra \hs  (ac)\Inf b = 
      (ca)\Inf b=0.  \nonumber \ee 
      Equivalent to (${\rm f}$) is the condition: \\[1mm] 
       \hspace*{1cm} 
         (${\rm f_1}$)  \;  $\{ab\}^{\perp\perp}\subset \{a\}^{\perp\perp}\cap\{b\}^{\perp\perp}\,$ for $\,0\leq a,b\in\cA$, \\see \cite[Proposition 3.5]{deP81}. 
\item An element $e\in \cA$ is called a \textit{multiplicative unit}, if  
      $a\cdot e = e\cdot a = a$ for all $a\in \cA$. It is uniquely defined. An algebra with 
      multiplicative unit is called {\it unitary}.
\item An element $a\in \cA$ is called \textit{nilpotent}, if there is $n\in\bN$ such that $a^n = 0$.
      The set of all nilpotent elements of $\cA$ is denoted by $N(\cA)$. If $\cA$ is an Archimedean $f$-algebra, then
      $N(\cA) = N_2(\cA) := \left\{ a \in\cA \colon a^2 =0\right\}$, see \cite[Proposition~10.2~(i)]{deP81}.
\item An $\ell$-algebra $\cA$ is called \textit{semiprime}, if the only nilpotent element in 
      $\cA$ is zero.
\end{itemize}
\rem\;{\bf 1}\;\label{f-und-d}
We collect here without proof the main properties of the introduced $\ell$-algebras and comment the relations between them. 
For the proofs we refer to \cite{AliBur85,BerHui90,Mey91,Zaa83}. 
Let $\cA$ be an arbitrary $\ell$-algebra. 
\begin{enumerate}
\item[(1)] It follows immediately from the definitions that each $f$-algebra is a $d$-algebra.
	The converse, in generally, is not true.
\item[(2)] If a $d$-algebra is semiprime or possesses a positive multiplicative unit, then it is an $f$-algebra.
\item[(3)] If in an $f$-algebra a multiplicative unit exists, then the latter is always positive.
\item[(4)] Even in an $f$-algebra the existence of a multiplicative unit is not guaranteed: 
	The vector lattice $c_0$ of all real zero sequences with the coordinatewise order and algebraic 
        operations is a semiprime Archimedean $f$-algebra without a multiplicative unit.
\item[(5)] An Archimedean $\ell$-algebra with a multiplicative unit $e>0$ is an $f$-algebra if and only 
        if $e$ is a weak order unit.
\item[(6)] Every Archimedean $f$-algebra is commutative and every unitary Archimedean $f$-algebra is semiprime.
\item[(7)] In an Archimedean commutative $d$-algebra
        the following frequently used formulas hold (see \cite[Proposition 1]{BeuHui84}, and \cite[Proposition 4]{Bou00})
	for the vector lattice ope\-ra\-tions with $p$-th powers of $a,b\in \cA_+$ for $p\in\bN_{\geq 1}$:
        \be\label{potenz} (a\Inf b)^p = a^p\Inf b^p \quad\mbox{and}\quad (a\Sup b)^p = a^p\Sup b^p.\ee
\item[(8)] In any $\ell$-algebra the condition $a\Inf b = 0 \;\Rightarrow\; ab =0$ and the condition $a^2 = |a|^2$ are equivalent (see \cite[Proposition 1.3]{BerHui90}). They hold in every $f$-algebra.
\end{enumerate} 
\par\medskip 
The following definitions are basic. 
\begin{Def}
Let $V$ be an Archimedean vector lattice. 
\begin{enumerate} 
\item\label{dfin} An element $\varphi\in V$ is called {\em finite},
	if there exists an element $z\in V_{+}$ such that the following condition holds:
	for any $x \in V$ there is a number $c_{x}>0$ satisfying the inequality 
	\[\left|x\right|\wedge n\left|\varphi\right|\leq c_{x}z\; \mbox{ for all }\; n \in \bN. \]
	The element $z$ is called a {\em $V$-majorant} or briefly a {\em majorant} of $\varphi$. 
\item\label{totfin}
	An element $\varphi$ of a vector lattice $V$ is called {\em totally finite},
	if it possesses a $V$-majorant which itself is a finite element. 
\item\label{selbstmdef}
	An element $\varphi\in V$ is called {\em self-majorizing}, 
	if $\left|\varphi\right|$ is a majorant of $\varphi$, i.e. for each element $x \in V$ there is a constant 
	$c_{x}>0$ such that 
	\be\label{f2} 
	\left|x\right|\wedge n\left|\varphi\right| \leq c_{x}\left|\varphi\right|\;\mbox{ for all }\; n \in \bN.
	\ee 
\end{enumerate}
\end{Def}
The sets of all finite elements and all totally finite elements in $V$ are denoted by $\Phi_{1}(V)$ and $\Phi_{2}(V)$, respectively. 
It is easy to see that $\Phi_{1}(V)$ and $\Phi_{2}(V)$ are order ideals in $V$ and $\Phi_{2}(V)\subseteq \Phi_{1}(V)$.
The set of all self-majorizing elements is denoted by $S(V)$, the set of positive self-majorizing elements by $S_{+}(V)$, i.e. $S_{+}(V)=S(V)\cap V_+$. It is clear that with $\varphi$ also $\bl \varphi \br$ is a self-majorizing element. The set $\Phi_{3}(V)=S_+(V)-S_+(V)$   
is an order ideal in $V$ and $\Phi_3(V)\subseteq\Phi_2(V)$. 
\par 
The main characterization of self-majorizing elements is contained in the following theorem. 
The proofs of the theorem and its corollary are provided as Theorem 1 and Corollary 3 in \cite{TeiWeb11}. 
\par
\begin{theo}[{\cite[Corollary~7.2]{LuxMoo67}}, and \cite{FelPor76}]\leavevmode\vspace*{1mm}
\label{T1}
For an element $\varphi$ of a vector lattice $V$ the following statements are equivalent.
\begin{enumerate}
	\item The element $\varphi$ is self-majorizing.
	\item The order ideal $I_{\varphi}$ generated in $V$ by $\varphi$ is the projection band $\left\{\varphi\right\}^{\perp\perp}$.
\end{enumerate}
\end{theo}
Theorem \ref{T1} yields the following corollary as an immediate consequence. 
\begin{cor}\label{C1} 
Let $V$ be a vector lattice. Then 
 \begin{enumerate}
	\item any order unit in $V$ is a self-majorizing element and
        \item if  $V$ possesses an order unit then 
              \mbox{$\Phi_{3}(V)=\Phi_{2}(V)=\Phi_{1}(V)=V$}.
 \end{enumerate}
\end{cor} 
\section[unitary l-algebras]{Finite elements in unitary $\ell$-algebras}\label{Sec3}
The first result shows that the multiplication with elements from the order ideal generated by the positive multiplicative unit preserves the finiteness with the same majorant. 
\begin{theo}\label{T2}
Let $\cA$ be an $\ell$-algebra with a positive  multiplicative unit $e$ and let $a$ be an arbitrary element of $I_e=\left\{a\in \cA \colon |a|\leq\lambda e \mbox{ for some } \lambda\in\bR_{\geq 0}\right\}$.
Then for $i=1,2$ there holds 
\[\varphi\in\Phi_i(\cA) \mbox{ with the majorant }  u \;\Longrightarrow \; 
			\varphi a,\, a\varphi \in\Phi_i(\cA) \mbox{ with the majorant } u.  \]
\end{theo}
\prf
Without loss of generality let $\varphi\geq 0$
(otherwise use $\varphi = \varphi^+-\varphi^-$).
It suffices to consider only $a\geq 0$, since by condition ($\ell_2$) there holds $|a\varphi| \leq |a|\varphi$.
For an $a\in I_e$,
there is a $\lambda\in\bR_{\geq0}$ such that $0\leq a\leq\lambda e$.
Due to the condition ($\ell_1$) we have for arbitrary $x\in\cA$ and all $n\in\bN$ 
the inequality 
\[ \bl x\br\Inf na\varphi \leq \bl x \br\Inf n\lambda e \varphi = \bl x\br \Inf n\lambda \varphi.\]
If now $\varphi$ is a finite element with a majorant $u$, then 
\[ \bl x\br\Inf na\varphi\leq c_x u \,\mbox{ for all } \,n\in\bN.\]
Therefore the product $a\varphi$ is also a finite element with the same majorant as $\varphi$.
Ana\-lo\-gous\-ly, the statement is proved for the product $\varphi a$.

If $\varphi$ is even totally finite, i.e. the majorant $u$ of $\varphi$ itself is a finite element, 
then the products $\varphi a$ and $a \varphi$ also have finite majorants, which shows that they are totally finite as well.
\qed

\medskip 
The same result can be proved without the positivity of the multiplicative unit, if $\cA$ is supposed to be a 
$d$-algebra. However, in contrast to the previous theorem, the majorant for the product changes and depends on the factor $a$. 
\begin{theo}\label{T3}
Let $\cA$ be a $d$-algebra with a (not necessarily positive) multiplicative unit and let $a\in \cA$ be an  arbitrary element.
Then for $i=1,2$ there holds
\[ \varphi\in\Phi_i(\cA) \;\Longrightarrow\; a \varphi,\hs \varphi a \in \Phi_i(\cA).
\]
In particular, $\Phi_i(\cA)$ is a $d$-subalgebra and a ring ideal in $\cA$. \\
If additionally $\cA$ is an $f$-algebra, then $\Phi_i(\cA)$ is even an $f$-subalgebra.
\end{theo}
\prf
Denote the multiplicative unit of $\cA$ by $e$ and
assume again $\varphi\geq 0$.
Let first $i=1$.
Due to $\varphi\in \Phi_1(\cA)$ there are a majorant $u\in \cA$ for $\varphi$ and, for each 
$x\in \cA_+$, a number $c_x\in\bR_{\geq 0}$ such that 
\be\label{f3} 
             \bl x\br\Inf n\varphi \hs\leq\hs c_x u \; \mbox{ for all }\; n\in\bN .
\ee
Since by condition (${\rm d_1}$) one has $|\varphi a| = |\varphi||a|$, the elements $\varphi |a|$ and $\varphi a$ 
are coincidently finite, so $a\geq 0$ may be assumed.
\par
Let $x\in \cA$ and $n\in\bN$ be arbitrary. Then $a\geq 0$ implies $a\Sup e\geq 0$ and 
by means of condition (${\rm d_3}$) one obtains from (\ref{f3})  
\[    (a\Sup e)\bl x\br \hs\Inf\hs n(a\Sup e)\varphi \hs\leq\hs c_x (a\Sup e)u.   \]
Since $a\leq a\Sup e$ and $\bl x\br= e\bl x\br \leq (a\Sup e) \bl x\br$ one has for every $n\in\bN$
\[ \bl x\br\hs\Inf\hs na\varphi \hs\leq\hs \bl x\br \hs\Inf\hs n(a\Sup e) \varphi \hs\leq\hs (a\Sup e)\bl x\br\hs\Inf\hs n(a\Sup e)\varphi \hs\leq\hs 
c_x (a\Sup e)u,\]
i.e. the element $a\varphi$ is  finite in $\cA$ with the majorant $(a\Sup e)u$.
Analogously it can be shown that the product $\varphi a$ is finite in $\cA$. 
\par
The set $\Phi_1(\cA)$ is an order ideal in $\cA$, in particular a vector sublattice. According to the first part of the proof the product of two finite elements is finite and so, the set $\Phi_1(\cA)$ is closed under the multiplication. The properties (d) or (f) are shared by the set $\Phi_1(\cA)$, if $\cA$ is a $d$- or an 
$f$-algebra, respectively. Therefore, 
$\Phi_1(\cA)$ is a $d$- or an $f$-subalgebra of $\cA$, respectively. \\
It is clear from the proof that $\Phi_1(\cA)$ is a ring ideal.
\par
For $i=2$ observe that $(a\Sup e) u$ is a majorant of the element $a\varphi$ as above, where $u$ as a majorant of the totally finite element $\varphi$ can be assumed to be a finite element. 
By what has been proved in the first part (case $i=1$)  the element $(a\Sup e) u$ is finite as well, which yields the totally finiteness of $a\varphi$ in $\cA$. 
The totally finiteness of the product $\varphi a$ is proved analogously. 
The remaining statements for $\Phi_2(\cA)$ are obtained analogously to the case  $i=1$.
\qed
\vspace*{3mm}\\
\rem\;{\bf 2}\;
If a majorant of $\varphi$ is $u$, then a majorant of $a\varphi, \, \varphi a$ is $(a\Sup e)u, \, u(a\Sup e)$, respectively. In particular, the idea of this proof cannot be used to obtain an analogous result for self-majorizing elements.
\par\medskip
If the multiplicative unit itself is a finite element then we get  
\begin{theo}\label{T4}
Let $\cA$ be a $d$-algebra with a multiplicative unit $e$. Let be $e\in \Phi_1(\cA)$.  Then 
\[ \Phi_1(\cA) \hs = \hs \Phi_2(\cA)\hs = \hs \cA.\]
If $\cA$ is an $f$-algebra and $e\in\Phi_1(\cA)$, then $e$ is even an order unit in $\cA$ and  
\[   \Phi_1(\cA)\hs = \hs \Phi_2(\cA)\hs = \hs \Phi_3(\cA)\hs = \hs \cA.  \]
\end{theo}
\prf
First consider the case of a $d$-algebra. Since $e$ is finite, by the previous theorem the products $a e$ and $e a$ are finite elements for all $a\in \cA$, i.e. $\cA\subseteq \Phi_1(\cA)$. So the equalities 
$\Phi_1(\cA)=\Phi_2(\cA)=\cA$ hold.  
\par
Consider the case of an $f$-algebra $\cA$. 
We mention first that the element $e$ is positive and a weak order unit [Remark 1(3) and 1(5)]. 
If $\cA$ has a weak order unit then, according to 
Corollary 2.5 from \cite{ChenWeb03-1}, the equalities $\Phi_1(\cA)\hs = \hs \Phi_2(\cA)\hs = \cA$ hold if and only if there exists an order unit in $\cA$. 
Since the equalities hold by what has been proved in the first part (here we use the fact that an $f$-algebra is also a $d$-algebra), the $f$-algebra $\cA$ has an order unit. 
From Corollary \ref{C1} it follows now that 
$\cA$ coincides also with the order ideal $\Phi_3(\cA)$ of all self-majorizing elements of $\cA$.  
\par
Now consider the weak order unit $e$ which, due to $\cA=\Phi_3(\cA)$,  is a self-majorizing element, and show that $e$ is a (strong) order unit. By Theorem \ref{T1} the order ideal generated in $\cA$ by $e$  is a projection band and coincides with $\{e\}^{\perp\perp}$. Since $e$ is a weak order unit one has $\cA=\{e\}^{\perp\perp}$. Consequently, $\cA=\{e\}^{\perp\perp}=\{a\in \cA\colon \exists \la>0 \mbox{ with } \bl a\br\leq \la e\}$, i.e. $e$ is an order unit in $\cA$. 
\qed
\par\medskip
\begin{theo}\label{T5}
Let $\cA$ be an $f$-algebra with a multiplicative unit $e$. 
Let there exist a submultiplicative Riesz norm on $\cA$, i.e. a Riesz norm which satisfies $\nl ab\nr\leq \nl a\nr \nl b\nr$ for all $a,b\in \cA_+$. Then
\begin{enumerate}
\item the multiplicative unit $e$ is an order unit and 
\item $\Phi_1(\cA)\hs=\hs\Phi_2(\cA)\hs=\hs\Phi_3(\cA)\hs=\hs\cA$. 
\end{enumerate}
\end{theo}
\prf\footnote{The main idea of the proof is due to W.A.J.~Luxemburg, cf. \cite[Theorem 15.5]{AliBur85}.}
1. 
We show that for each $a\in\cA, \ 0\neq a$ there exists a $\lambda\in \bR_{\geq 0}$ such that $-\lambda e \leq a\leq \lambda e$.
Let first $a\in \cA_+$. Further on we use the obvious decomposition \mbox{$a-\lambda e = (a-\lambda e)^+ - (a-\lambda e)^-$}, which holds for any $\lambda\in \bR_{\geq 0}$.
Without loss of ge\-ne\-ra\-li\-ty $(a-\lambda e)^+>0$ can be assumed. Indeed, $(a-\lambda e)^+ = 0$ for some $\lambda>0$ leads to $a-\lambda e = -(a-\lambda e)^-\leq 0$, and so to $0\leq a\leq \lambda e$. 
\par
Now consider the element $(a-\lambda e)(a-\lambda e)^+$.
Due to Remark 1(8) the product of the two positive disjoint elements $(a-\lambda e)^-$ and $(a-\lambda e)^+$ vanishes, and by taking into account the condition ($\ell$) we obtain the inequality
\[(a-\lambda e)(a-\lambda e)^+ = (a-\lambda e)^+(a-\lambda e)^+ \hs-\hs (a-\lambda e)^-(a-\lambda e)^+=
	 \big((a-\lambda e)^+\big)^2 \geq 0.\]
We conclude $a(a-\lambda e)^+ - \lambda e (a-\lambda e)^+ \geq 0$, and so 
\[a(a-\lambda e)^+ \geq \lambda (a-\lambda e)^+ >0 \;\mbox{ for } \la>0. \]
Due to the norm being submultiplicative and Riesz we obtain
\[\nl a \nr \nl (a-\lambda e)^+ \nr \geq \nl a (a-\lambda e)^+ \nr \geq \lambda \nl (a-\lambda e)^+\nr >0, \] 
and therefore $\lambda \leq \nl a\nr$.
\par
Altogether, as we have seen, the assumption $(a-\lambda e)^+>0$ leads to $\lambda \leq \nl a\nr$.
Therefore $\lambda > \nl a\nr$ yields $(a-\lambda e)^+=0$, and so $a-\lambda e = -(a-\lambda e)^-\leq 0$, and again $0\leq a\leq \lambda e$ as above.
\par
Now let $a\in\cA$ be an arbitrary element. In view of $ \pm a\leq |a|$ we obtain the claimed result.
\par
2. The fact that all elements in $\cA$ are finite, totally finite and even self-majorizing follows from Corollary \ref{C1} by taking into account that $e$ is an order unit in $\cA$.
\qed
\par\bigskip
A linear operator $T$ on an Archimedean vector lattice $V$ is called {\it band preserving} if 
$T(B)\subseteq B$ for each band $B$ in $V$. A band preserving operator which is order bounded is called an {\it orthomorphism}. 
\par
It is well known that the collection \Orth($V$) of all orthomorphisms on an Archi\-me\-dean vector lattice $V$ is an $f$-algebra with the identity as a weak order unit. Moreover, any $f$-algebra $\cA$ with a multiplicative unit $e$ is algebraic and lattice isomorphic to \Orth$(\cA)$, where the image of $e$ is the identity in Orth$(\cA)$ (\cite[Theorems 3.1.10 and 3.1.13]{Mey91}).   
\begin{cor}\label{C2}
Let $\cA$ be a unitary $f$-algebra. 
Let there exist a submultiplicative Riesz norm on $\cA$. Then
\begin{enumerate}
\item the identity operator $I$ is an order unit in ${\rm Orth}(\cA)$ and 
\item ${\rm Orth}(\cA) \hs=\hs\Phi_i\big({\rm Orth}(\cA)\big), \hs\hs\hs i=1,2,3$. 
\end{enumerate}
\end{cor}
A similar result holds if there is some norm on the algebra $\cA$ which turns it into a Banach lattice. 
\begin{theo}\label{T6}
Let $\cA$ be an $f$-algebra with a multiplicative unit $e$. 
Let there exist a norm on $\cA$ such that $\cA$ becomes a Banach lattice. Then 
\begin{enumerate}
\item the multiplicative unit $e$ is an order unit and 
\item $\Phi_1(\cA)\hs=\hs\Phi_2(\cA)\hs=\hs\Phi_3(\cA)\hs=\hs\cA$. 
\end{enumerate}
\end{theo}
\prf 
Since the $f$-algebra $\cA$ and Orth$(\cA)$ are algebraic and lattice isomorphic such that the image of $e$ under the isomorphism is $I\in\mbox{Orth}(\cA)$, then by Wickstead's Theorem (\cite[Theorem 15.5]{AliBur85}) the identity operator $I$ is an order unit in Orth($\cA$) and so, $e$ is an order unit in $\cA$. 
By virtue of Corollary \ref{C1} all elements in Orth($\cA$), and consequently in $\cA$, are finite, 
totally finite and even self-majorizing.  
\qed
\medskip

For the $f$-algebra of all orthomorphisms on a vector lattice from the Theorems \ref{T3} and \ref{T6} we get the following properties which we formulate as
\begin{cor}\label{C3}
\leavevmode
\begin{enumerate}
\item Let $V$ be a vector lattice. If $S\in \Phi_i\big({\rm Orth}(V)\big)$ for $i=1,2$ and $T\in {\rm Orth}(V)$  
then also $S\circ T\in\Phi_i({\rm Orth}(V))$. In particular, $\Phi_i\big({\rm Orth}(V)\big)$ is an $f$-sub\-al\-geb\-ra and a ring ideal. 
\item Let $E$ be a Banach lattice. Then ${\rm Orth}(E)$ is an $f$-algebra and under the order unit  norm $\nl T\nr_I=\inf\{\la>0\colon \bl T\br\leq \la I\}$ also an AM-space with order unit. In this case
\[{\rm Orth}(E)\hs=\hs\Phi_i\big({\rm Orth}(E)\big)\quad \textnormal{for }\; i=1,2,3.\] 
\end{enumerate}
\end{cor}
\par
The last results throw some light also on the relations between finiteness and invertibility of elements in $f$-algebras.
\par\smallskip
\hspace*{1cm}\\
\bsp{\bf 1}\;\label{ex1} 
Consider the vector lattice $C_b(\bR)$ of all bounded real-valued continuous functions on $\bR$ equipped with the pointwise algebraic operations and partial order. 
Then $C_b(\bR)$ turns out to be an Archimedean $f$-algebra.  
It is a Banach lattice if the norm is defined by $\nl f\nr_{\infty}=\sup\limits_{x\in\bR}|f(x)|$ for 
$f\in C_b(\bR)$. 
Since there exist (many) order units in $C_b(\bR)$ all elements are finite\footnote{
The Banach algebras $C_b(\bR)$ and $C(\beta\bR)$ are lattice isomorphic, where $\beta\bR$ denotes the 
Stone-\v{C}ech compactification of $\bR$. So, all elements in $C_b(\bR)$ like in $C(\beta\bR)$ are finite.}. 
 Observe that any function $f\in C_b(\bR)$ with $\inf\limits_{x\in \bR}|f(x)|>0$ is invertible. Of course, there are non-invertible finite elements as well, e.g. functions with compact support.  
\par
Let $\cA$ be a $d$-algebra with a multiplicative unit $e$. If there exists at least one non-zero finite element which is invertible in $\cA$, then immediately all elements of $\cA$ are finite, i.e. $\cA=\Phi_1(\cA)$. Indeed, Theorem \ref{T3} guarantees that the finiteness and the invertibility of an element $\varphi$ 
imply $e=\varphi^{-1}\varphi$ to be a finite element in $\cA$. Then by Theorem \ref{T4} all elements of $\cA$ are finite.  
\par
The $f$-algebra $C(\bR)$ of all continuous functions on $\bR$ contains a multiplicative unit (the function \ou). However, in contrast to $C_b(\bR)$, there is no order unit.
There is also no norm on 
$C(\bR)$ that makes it a Banach lattice. Otherwise, by Theorem \ref{T6}, there would be an order unit. 
 It is clear that the element \ou \, is not finite\footnote{since $\Phi_1(C(\bR))={\mathcal K}(\bR)$, the vector lattice of all functions with compact support.} in $C(\bR)$. By what has been mentioned above no finite element can be invertible. Consequently, there exist $f$-algebras $\cA$ with multiplicative units such that $\Phi_1(\cA)\neq \{0\}$ and no finite element is invertible.   
\section[non-unitary f-algebras]{Finite elements in non-unitary $f$-algebras}\label{Sec4}
In this section we consider the case of $f$-algebras which do not possess any multiplicative unit. 
\par
\begin{Def}\label{q-close} 
Let $\cA$ be an $f$-algebra. 
\begin{enumerate}
\item $\cA$ is said to be {\rm square-root closed} if for any $a\in \cA_+$ there exists $b\in\cA$ such that 
$b^2=a$, i.e. for every such element $a$ there exists its square root. 
\item $\cA$ is said to have the {\rm factorization property} if for every $a\in \cA$ there exist two elements $b,c\in \cA$ such that $a=bc$. 
\item  We say that $\cA$ has the {\rm weak factorization property} if for every $a\in \cA$ there exist two elements $b,c\in \cA$ such that $a\leq bc$. 
\end{enumerate}
\end{Def}
In \cite[Theorem 4.6]{BeuHuideP83} the first two properties were proved to be equivalent in uniformly complete $f$-algebras. The fact that the property 3. is weaker than 2. is demonstrated by the next example.  
\par\smallskip
\hspace*{1cm}\\
\bsp{\bf 2}\;\label{ex2} 
For the vector lattice $\cA := \left\{f\in C[-1,1] : f(0) =0\right\}$ let the multiplication for all $f,g\in\cA$ be defined by
\[(f\cdot g)(t) :=\begin{cases}
  f(t)g(t), 	& t\in[0,1],\\
  f(-t)g(-t),	& t\in[-1,0).
\end{cases}\]
Products in $\cA$ are precisely the axisymmetric functions, which vanish at $0$.  
Observe that $\cA$ is an $f$-algebra, which is not semiprime.
We will show that $\cA$ is uniformly complete and has the weak factorization property. However, the factorization property does not hold for $\cA$.
\par
To see that $\cA$ is uniformly complete, notice that $\cA$ is the kernel $\delta_0^{-1}(0)$ of the continuous functional $\delta_0$ defined on the Banach lattice $C[-1,1]$ by $\delta_0(f)=f(0)$.
\par
The $f$-algebra $\cA$ obviously does not have the factorization property, since an arbitrary $g\in\cA$, which is not axisymmetric, cannot be written as a product of two elements of $\cA$. Since in uniformly complete $f$-algebras the factorization property is equivalent to the square-root closedness, the latter does not hold in $\cA$ either.
However, $\cA$ has the weak factorization property. Indeed, let $g\in\cA$ be an arbitrary element. Define
\[\hat{g}(t):=\mathop{\max}\limits_{t\in\left[-1,1\right]}\left\{|g(t)|,|g(-t)|\right\} \quad\text{ and }\quad \tilde{g}(t):=\sqrt{\hat{g}(t)}.\]
Then $\hat{g},\tilde{g}\in\cA$ and there holds the inequality $g\leq \hat{g} = \tilde{g}^2$.
\par\smallskip 
\hspace*{1cm}\\
\bsp{\bf 3}\;\label{ex3} 
This example shows that, in general, the weak factorization property does not hold in $f$-algebras.
To that end, consider finite partitions $\tau$ of the set $[0,\infty)$, i.e. $\tau=\{I_0,\ldots,I_n\}$ such that 
\[\bigcup\limits_{k=0}^n I_k=[0,\infty),\quad I_k\cap I_j=\emptyset\,\,\,\,\textnormal{ for } k\neq j\]
and $I_k$ is a subinterval of $[0,\infty)$ for any $k$.  
Let $\cP$ be the set of all polynomials $p$ vanishing at the point $t=0$.   
Consider now the collection $\fP:=\fP\big([0,\infty)\big)$ of all continuous functions 
on $[0,\infty)$ for each of which there exists a partition $\tau$ such that $f\big|_{I_k}=p_k$ 
with $p_k\in \cP$ for any $I_k\in \tau$. 
The algebraic operations, the multiplication 
and the partial order are introduced in $\fP$ pointwise. Then $\fP$ is an Archimedean
 $\ell$-algebra. Moreover, it 
is easy to see that the disjointness\footnote{The supports of two disjoint continuous functions on $[0,\infty)$ intersect at most in one point.} of two functions $f, g\in \fP$ is preserved also after the 
multiplication of one of them by a positive function $h\in \fP$. 
Therefore $\fP$ is an $f$-algebra. Since only the zero-element of $\fP$ can satisfy the equation $f^2=0$, 
the $f$-algebra is semiprime. 
Observe that the restriction on $[0,\infty)$ of a polynomial $p$ of arbitrary degree with $p(0)=0$ belongs 
to $\fP$ but the function \ou$\big|_{[0,\infty)}$ does not. It follows that $\fP$ does not contain neither 
an order unit nor a multiplicative unit. 
\par
The $f$-algebra $\fP$ does not have  even the weak factorization property, since the function $f(t)=t$ can not be estimated by 
a product  of two functions. 
Indeed, $f\leq pq$ implies that both polynomials $p,q$ take on positive values for all $t>0$ and ${\rm deg}(pq)\geq 2$.  
Since $pq(0)=0$ the graphs of $f$ and $pq$ intersect in some point. Let $t_0$ be 
the smallest number with $0<t_0$ and $f(t_0)=pq(t_0)$. There is an interval $I_k$ 
of a partition for $pq\in \fP$ such that $t_0\in I_k$ and $f(t)> pq(t)$ for $t\in (0,t_0)$. 
\par\bigskip
The $p$-fold product is used to define the $p$-th root of an element in an $\ell$-algebra: for $g\in \cA$ an element $\tilde{g}\in \cA$ is called a {\it ${p}$-th root} of $g$,
if $\tilde{g}^p = g$. If a $p$-th root of $g$ exists and is uniquely defined then we write $\tilde{g} = g^{\frac{1}{p}}$ and call $\tilde{g}$ \textit{the} $p$-th root of $g$. 
For details we refer to \cite{BeuHui84,Bou00,Zaa83} where, in particular, the following results can be found.
\\[2mm]
\rem\;{\bf 3}\;\label{root} 
\vspace{-2mm}
{\begin{enumerate}
\item[(1)]{(Existence and uniqueness of the root).}
	Let $\cA$ be an Archimedean uniformly complete $f$-algebra and 
	$p\in \bN_{\geq 2}$.
	Then (see \cite[Theorem 3]{Bou00}) there exists a positive $p$-th root for any $p$-fold product\footnote{For the product of $p$ elements $g_1,\dots, g_p$ in $\cA$ we will use the notation $g_1 \cdots g_p := g_1\cdot \ldots \cdot g_p$.} of positive 
        elements of $\cA$, i.e.
\be\label{f4} 
  g_1,\dots, g_p \in \cA_+ \; \Longrightarrow\; (g_1 \cdots g_p)^{\frac{1}{p}}
							\mbox{ exists in } \cA_+.
\ee
The root is uniquely defined if the algebra $\cA$ is semiprime. 
\item[(2)]{(Monotonicity of the root).} In every $\ell$-algebra $\cA$ for $p\in \bN_{\geq 2}$ and $a,b \in \cA_+$ due to the property $(\ell_1)$ one has
        \[a\leq b \;\Longrightarrow \;a^p\leq b^p.\] 
If $\cA$ is a semiprime $f$-algebra, then (see \cite[Theorem 142.3]{Zaa83} and \cite[Proposition 2.(iii)]{BeuHui84}) the root is monotone, i.e.
        \[a\leq b \;\Longleftrightarrow \;a^p\leq b^p.\] 
\end{enumerate}
\par
\begin{theo}\label{T7}
Let $\cA$ be a uniformly complete $f$-algebra with the weak factorization property, $p\in\bN_{\geq2}$ and $i=1,2$. If $\varphi_1,\dots,\varphi_p \in \Phi_i(\cA)$ with majorants $u_1,\ldots,u_p$, respectively,
then $\varphi_1\cdots\varphi_p  \in\Phi_i(\cA)$ with a majorant $(u_1\Sup \ldots \Sup u_p)^p$.
In particular, $\Phi_i(\cA)$ is an $f$-subalgebra of $\cA$.
\end{theo}
\prf
Let $i=1$. First we prove the claim for the $p$-fold power $\varphi^p$ of a finite element $\varphi\in\cA$.
Let $\varphi$ be a finite element in $\cA$ with a majorant $u\in \cA_+$. 
Without loss of generality, $\varphi$ can be assumed to be positive, otherwise use $|\varphi^p|=|\varphi|^p$, which holds due to property (d$_1$).
For an arbitrary $a\in \cA$ the weak factorization property of $\cA$ yields the existence of $p$ elements $a_1,\dots, a_p \in \cA$ with $|a|\leq a_1\cdots a_p$.
Again by property ($\ell_2$) and Remark 3(1) there follow the positivity of the elements $a_1,\dots,a_p$ and the existence of a root $(a_1\cdots a_p)^\frac{1}{p}$ in $\cA$.
Using the formula (\ref{potenz}) and the finiteness of $\varphi$ we obtain that there is a constant 
$c_{\sqrt[p]{a_1\cdots a_p}}\geq 0$ such that for all $n\in\bN$
\[ |a| \wedge n\varphi^p \hs\leq\hs (a_1\cdots a_p) \wedge n\varphi^p
	\hs=\hs \big((a_1\cdots a_p)^\frac{1}{p} \wedge n \varphi\big)^p \hs\leq\hs c^p_{\sqrt[p]{a_1\cdots a_p}}\hs u^p,\]
where the last inequality follows from Remark 3(2).
Therefore the $p$-th power $\varphi^p$ of a finite element $\varphi\in\cA$ is also finite.
\par
Now let $\varphi_1,\dots,\varphi_p$ be arbitrary finite elements in $\cA$ with majorants 
$u_1,\ldots, u_p$, respectively.
The modulus of the product $\varphi_1\cdots\varphi_p$ can be estimated by
\[ |\varphi_1\cdots\varphi_p| \hs\leq\hs |\varphi_1|\cdots|\varphi_p| \hs\leq\hs \underbrace{(|\varphi_1| \vee \dots \vee|\varphi_p|)\cdots(|\varphi_1| \vee\dots\vee |\varphi_p|)}_{p \textnormal{ times}} = (|\varphi_1| \vee\dots\vee |\varphi_p|)^p.\]
Since a majorant of the supremum $|\varphi_1| \vee\dots\vee |\varphi_p|$ is given by $u_1\vee \ldots\vee u_p$, by the first part the $p$-fold product $(|\varphi_1| \vee\dots\vee |\varphi_p|)^p$ is finite as well with the majorant $(u_1\vee \ldots\vee u_p)^p$.
\par\smallskip
Let $i=2$. Then the majorants $u_1,\ldots,u_p$ can be assumed to belong to $\Phi_1(\cA)$ and 
the element $(u_1\vee\ldots\vee u_p)^p$ is finite due to what has been proved in the case $i=1$. 
\qed
\par\medskip
The last theorem has been proved under stronger conditions than Theorem \ref{T9}, where we drop the uniformly completeness and the factorization property of the $f$-algebra. However, in the proof of Theorem \ref{T9} the majorants are not given explicitly and so, in contrast to Theorem \ref{T7}, the fate of totally finite elements remains unknown there.
\par\medskip
For the next theorem notice that the Example 2 shows that the weak factorization property does not imply semiprimitivity, even under the additional condition of uniformly completeness. By Example 3 the converse implication is also not true. However, it is not known if uniformly completeness together with semiprimitivity imply the weak factorization property.
\par\medskip
\begin{theo}\label{T8}
Let $\cA$ be a semiprime uniformly complete $f$-algebra with the weak factorization property and $p\in\bN_{\geq2}$.
If for $\varphi\in \Phi_1(\cA)$ there exists the root $\varphi^{\frac{1}{p}}$ in $\cA$, then 
 $\varphi^{\frac{1}{p}}\in\Phi_1(\cA)$.
\end{theo}
\prf 
First consider $0<\varphi\in \Phi_1(\cA)$ with a majorant $u\in \cA$ for which there exists the root 
$\varphi^{\frac{1}{p}}$ in $\cA_+$.
Let $a\in \cA_+$ be an arbitrary element. 
According to the formula (\ref{potenz}) and by using the finiteness of $\varphi$ we get 
\[ (a\Inf n\varphi^{\frac{1}{p}})^p \hs=\hs a^p \Inf n^p \varphi \hs\leq\hs c_{a^p} u\]
for some constant $c_{a^p}\geq 0$ and all $n\in \bN$.
Due to the weak factorization property of $\cA$ there are $p$ elements  $u_1,\dots, u_p \in \cA$ such 
that $u\leq u_1\cdots u_p$.
Therefore the above inequality can be continued as follows 
\[ (a\Inf n\varphi^{\frac{1}{p}})^p \hs\leq\hs c_{a^p} u \hs\leq\hs c_{a^p} u_1\cdots u_p.\]
Due to condition ($\ell_2$) there hold the relations   
$0\leq u\leq u_1\cdots u_p = |u_1\cdots u_p|\leq |u_1|\cdots|u_p|$. 
Without loss of generality we may replace $u_i$ by $\bl u_i\br$ and therefore assume that $u_i\geq 0$ for $i=1,\ldots, p$.  
According to  Remark 3(1) there exists the root $(u_1\cdots u_p)^\frac{1}{p}$ in $\cA$
and we obtain
\[ (a\Inf n\varphi^{\frac{1}{p}})^p \hs\leq\hs c_{a^p} u_1\cdots u_p = c_{a^p}((u_1\cdots u_p)^\frac{1}{p})^p.\]
The monotonicity of the root allows us to extract the $p$-th root on both sides, which yields  
\[ a\Inf n\varphi^{\frac{1}{p}} \hs\leq\hs \sqrt[p]{c_{a^p}}(u_1\cdots u_p)^\frac{1}{p},\] 
and that shows that the element $\varphi^{\frac{1}{p}}$ is finite in $\cA$ with the majorant $(u_1\cdots u_p)^\frac{1}{p}$.
\par 
Now let $\varphi\in\Phi_1(\cA)$ be arbitrary. If $\varphi$ possesses a root 
$\varphi^{\frac{1}{p}}$ then by condition ($\rm{d}_1$)
\[|\varphi|=|\varphi^{\frac{1}{p}}\cdots\varphi^{\frac{1}{p}}| = |\varphi^{\frac{1}{p}}|\cdots|\varphi^{\frac{1}{p}}| \qquad (p {\textnormal{ times}})\] 
implies that $|\varphi^{\frac{1}{p}}|$ is the $p$-th root of 
$|\varphi|$, i.e. $|\varphi^{\frac{1}{p}}|=|\varphi|^{\frac{1}{p}}$.
Together with $\varphi$  the element $|\varphi|$ is also finite in $\cA$ and so, according to the first part of the proof, the element $|\varphi^{\frac{1}{p}}|$ is finite and, therefore the finiteness of 
$\varphi^{\frac{1}{p}}$ is obtained.
\qed
\par\medskip
In analogy to the above theorem we obtain the next result.
\begin{cor}\label{C4}
Let $\cA$ be a semiprime $f$-algebra and $p\in\bN_{\geq2}$.
If for $\varphi\in \Phi_3(\cA)$ there exists the root $\varphi^{\frac{1}{p}}$ in $\cA$, then 
 $\varphi^{\frac{1}{p}}\in\Phi_3(\cA)$.
\end{cor}
\prf
First consider $0<\varphi\in \Phi_3(\cA)$ for which there exists the root $\varphi^{\frac{1}{p}}$ in $\cA_+$.
According to the formula (\ref{potenz}) and since $\varphi$ is self-majorizing we get 
\[ (a\Inf n\varphi^{\frac{1}{p}})^p \hs=\hs a^p \Inf n^p \varphi \hs\leq\hs c_{a^p} \varphi\]
for some constant $c_{a^p}\geq 0$ and all $n\in \bN$.
The monotonicity of the root allows us to extract the $p$-th root on both sides, which yields 
\[ a\Inf n\varphi^{\frac{1}{p}} \hs\leq\hs \sqrt[p]{c_{a^p}}\varphi^\frac{1}{p},\] 
and that shows that the element $\varphi^{\frac{1}{p}}$ is self-majorizing.
\par
Now let $\varphi\in\Phi_3(\cA)$ be arbitrary such that $\varphi^{\frac{1}{p}}$ exists. The application of the identity  $|\varphi^{\frac{1}{p}}|=|\varphi|^{\frac{1}{p}}$ analogously to the proof of the previous theorem ensures that $\varphi^{\frac{1}{p}}$ is self-majorizing.
\qed
\par\medskip
For the next result, which is similar to Theorem \ref{T3}, we use the characterization of $f$-algebras given by the condition (${\rm f_1}$).
\begin{theo}\label{T9}
Let $\cA$ be an $f$-algebra,  $\varphi\in\Phi_1(\cA)$ and $a\in\cA$.
Then $a\varphi\in\Phi_1(\cA)$. In particular, $\Phi_1(\cA)$ is an $f$-subalgebra and a ring ideal.
\end{theo}
\prf
Let first $a\in\cA_+$ and $\varphi\in\Phi_1(\cA)$, $\varphi\geq 0$.
Using the condition ($\rm f_1$)  we obtain 
$\left\{a \varphi \right\}^{\perp\perp} \hs\subseteq\hs \left\{a\right\}^{\perp\perp} \cap \left\{\varphi\right\}^{\perp\perp}$.
By \cite[Theorem 2.4]{ChenWeb03-1} for a finite element we have $\left\{\varphi\right\}^{\perp\perp} \subseteq \Phi_1(\cA)$ and so
\[\left\{a \varphi \right\}^{\perp\perp} \hs\subseteq\hs \left\{a\right\}^{\perp\perp} \cap \left\{\varphi\right\}^{\perp\perp} \hs\subseteq\hs \Phi_1(\cA).\]
In particular, the product $a \varphi$ is finite in $\cA$.
\par
Now let $a\in\cA$ be arbitrary and $\varphi$ positive. The first part of the proof yields $a^+\varphi,a^-\varphi \in\Phi_1(\cA)$ and so we obtain the finiteness of  $a\varphi=a^+\varphi-a^-\varphi$ in $\cA$.
Finally, assume $\varphi$ to be arbitrary. Since $\Phi_1(\cA)$ is an order ideal, we obtain the finiteness of $\varphi^+$ and $\varphi^-$ and therefore also the finiteness of $a \varphi = a\varphi^+ - a\varphi^-$ in $\cA$.
\qed
\par\medskip
Notice that the product $a_1\cdots a_p$ belongs to $\Phi_1(\cA)$ if at least one of the elements $a_1,\ldots,a_p\in\cA$ belongs to $\Phi_1(\cA)$.
\par\medskip
The next theorem generalizes Theorem \ref{T6}, since, as was already mentioned in Remark 1(6), a unitary $f$-algebra $\cA$ is automatically semiprime.  
For its proof we need the following result, which we obtain by resuming and restricting \cite[Theorem 12.3.8.]{BKW77}. 
\par
First we introduce the following notation. Let $\cA$ be an $\ell$-algebra and $c\in \cA$. Denote by $_c\pi$ and $\pi_c$ the left and right multiplications by $c$, respectively, i.e. $_c\pi, \, \pi_c\colon \cA\to \cA,$ defined by 
\[  _c\pi(a)=c\,a \; \mbox{ and }\; \pi_c(a)=a\,c\; \mbox{ for all }\;a\in \cA.  \]   
It is clear that every multiplication operator $_c\pi, \pi_c$ is order bounded.
If $\cA$ additionally satisfies the condition (f), then for $c\geq 0$ the operators $_c\pi$ and $\pi_c$ are band preserving (and hence orthomorphisms), since then one has 
$\pi_c(a)\Inf b\hs=\hs_c\pi(a)\Inf b=0$ whenever $a\Inf b=0$ (see \cite[Theorem 8.2]{AliBur85}). 
\par
Notice that the map $h: a\mapsto \pi_a$ from an $f$-algebra $\cA$ into ${\rm Orth}(\cA)$ is   
a homomorphism. Indeed the condition (d$_2$) implies 
\[
\pi_{a\wedge b}(c) \hs=\hs (a\wedge b) c \hs=\hs ac \wedge bc \hs=\hs \pi_{a}(c) \wedge \pi_{b}(c) \hs=\hs (\pi_{a}\wedge \pi_{b})(c)
\]
and thus $h(a\wedge b) = h(a) \wedge h(b)$. The other properties of $h$ follow analogously.
\begin{Lem}\label{L1}
For an Archimedean $f$-algebra $\cA$ the following conditions are equivalent:
\begin{enumerate}
\item The algebra $\cA$ is semiprime.
\item The map $h$ is an injective homomorphism from $\cA$ into ${\rm Orth}(\cA)$.
      In particular, $\cA$ is embeddable as an $f$-subalgebra into the Archimedean 
      unitary $f$-algebra ${\rm Orth}(\cA)$.
\end{enumerate}
\end{Lem}
\prf
$\Rightarrow$: 
Since $\cA$ is semiprime, one has $\pi_a\neq 0$ for all 
$a\in\cA, \, 0\neq a$.
Therefore $\text{ker}(h)=\{0\}$, i.e. $h$ is injective.
\par
$\Leftarrow$:
Since $\cA$ is embeddable into ${\rm Orth}(\cA)$ by means of $h$, we can identify $\cA$ with a sublattice of ${\rm Orth}(\cA)$. Let $a$ be a nilpotent element in $\cA$. Then the element $a$ is also nilpotent in ${\rm Orth}(\cA)$. But the unitary $f$-algebra ${\rm Orth}(\cA)$ is
semiprime, i.e. $a$ is the zero element in ${\rm Orth}(\cA)$ and also in $\cA$.
\qed
\vspace*{3mm}\\
\rem\;{\bf 4}\;
Let $\cA$ be a semiprime $f$-algebra. Then
\[\pi_\varphi \in\Phi_3\big({\rm Orth}(\cA)\big) \hs\Longrightarrow\hs \varphi\in\Phi_3(\cA).\]
Indeed, by Lemma \ref{L1} the $f$-algebra $\cA$ is a sublattice of $\Orth(\cA)$, and so for each $x\in\cA$ we obtain
\[|x| \Inf n|\varphi| \hs=\hs |\pi_x| \Inf n|\pi_{\varphi}| \hs\leq\hs c_{|\pi_x|}|\pi_{\varphi}| \hs=\hs c_{|\pi_x|}\bl \varphi\br \]
for any $n\in \bN$ and some constant $c_{|\pi_x|}\in\bR_{\geq 0}$. 
Notice that the same statement for finite and totally finite elements, in general, is not true, since in these cases the majorants might not belong to $\cA$.
\par
The inverse implication, in general, is not true because for $\varphi\in \Phi_3(\cA)$ the element $\pi_\varphi\in {\rm Orth}(\cA)$ may not be a majorant 
for itself. Indeed, if $\cA$ does not possess a multiplicative unit then for  $x\in {\rm Orth}(\cA)\setminus \cA$ a corresponding constant $c_x$ might not exist.  
\par\medskip
\begin{theo}\label{T10} 
Let $\cA$ be a semiprime $f$-algebra and let there exist a norm on $\cA$, under which $\cA$ 
is a Banach lattice. 
Then 
\[\Phi_1(\cA) \hs=\hs \Phi_2(\cA) \hs=\hs \Phi_3(\cA) \hs=\hs \cA.\]
\end{theo}
\prf 
Since $\cA$ is semiprime, according to Lemma \ref{L1} 
the $f$-algebra $\cA$ can be embedded as a subalgebra into $\Orth(\cA)$. 
We write $\cA\subseteq \Orth(\cA)$ after identifying $\cA$ with its image $h(\cA)$ in $\Orth(\cA)$. 
According to \cite[Theorem 15.5]{AliBur85} the identity $I$ is an order unit in $\Orth(\cA)$. 
By Corollary \ref{C1} we obtain
\[ \Phi_3\big(\Orth(\cA)\big) \hs=\hs \Orth(\cA) \hs\supseteq\hs \cA.\]
It follows for two arbitrary elements $a,\varphi\in \cA_+$ that 
\[a \Inf n\varphi \hs\leq\hs c_a \varphi,\]
i.e. all positive elements in $\cA$ are self-majorizing. Therefore, as the cone $\cA_+$ is reproducing in $\cA$,  each element $x\in \cA$ is self-majorizing, and we get  
$\Phi_1(\cA) =\Phi_2(\cA) =\Phi_3(\cA)=\cA$.
\qed
\section[Product algebras]{Finite elements in product algebras}\label{Sec5}
Let $\cA$ be an $\ell$-algebra and $p\in \bN_{\geq 2}$.
The following construction is well-known. For details the reader is referred to \cite{BeuHui84,Bou00,Bou03,LuxZaa71}.
By 
\[ \Pi_p(\cA) := \left\{ g_1\cdots g_p \colon g_i \in \cA \mbox{ for } i=1,\dots,p \right\}\subseteq \cA\]
we denote the set of all $p$-fold products in $\cA$.
Clearly, $\Pi_p(\cA)\subseteq \cA$. In general, this inclusion is proper, e.g. define $\cA$ to be as in Example 2.
Even if $\cA$ is a semiprime uniformly complete $f$-algebra, then, in general, still $\Pi_p(\cA) \neq \cA$, see e.g. \cite{BeuHuideP83}, page 136, where an example for a semiprime uniformly complete and not square-root closed $f$-algebra is provided.
\par
If the set $\Pi_p(\cA)$, equipped with the order 
and algebraic operations induced from $\cA$, turns out to be an algebra, 
then it is called the 
\textit{product algebra of order $p$ of $\cA$.}
Denote by 
\[ \Sigma_p(\cA) := \left\{ g^p \colon g \in \cA_+ \right\}\]
the set of all $p$-fold powers of positive elements of $\cA$.
\par\smallskip
For completeness we provide without proofs some important properties of $\Pi_p(\cA)$. 
\begin{itemize}
\item If $\Pi_p(\cA)$ is a vector space, it may fail to be a vector lattice, in general.
	In case of $p=2$ there is a counterexample of an $\ell$-algebra $\cA$, which shows that the vector space $\Pi_2(\cA)$ is not a vector lattice under the order induced from $\cA$ (see \cite{Bou00}, Example 1).
\end{itemize}
Let $\cA$ be a uniformly complete $f$-algebra and $p\in\bN_{\geq 2}$.
\begin{itemize}
\item The set $\Pi_p(\cA)$ is a semiprime uniformly complete $f$-subalgebra of $\cA$ (see \cite[Corollary 5.3(iv)]{Bou03} \cite[Corollary 3]{Bou00} and \cite[Corollary 4]{BBT03}).
\item The set $\Pi_p(\cA)$ is a vector lattice under the ordering inherited from $\cA$, where
	\begin{equation}\label{g1} 
		\Pi_p^+(\cA)=\Sigma_p(\cA). 
	\end{equation}
	Additionally, for the supremum $\Sup_p$ and infimum $\Inf_p$ in $\Pi_p(\cA)$ the following formulas hold:
	\begin{equation}\label{f5} 
		f^p\Inf_p g^p = (f\Inf g)^p \quad\text{ and }\quad f^p\Sup_p g^p = (f\Sup g)^p \quad\text{ for }\quad f,g\in\cA_+.
	\end{equation}
\item If $f_1,\dots,f_p \in \cA$ are arbitrary elements, then for the modulus of the product $f_1\cdots f_p$  in $\Pi_p(\cA)$ the following formula is true  
	\be\label{f5a}   |f_1\cdots f_p|_p = |f_1|\cdots|f_p|. \ee 
	(see \cite[Proposition 1]{BeuHui84} and \cite[Corollary 5.3(i), (iv)]{Bou03}).
	\end{itemize}
Altogether we obtain for a uniformly complete $f$-algebra $\cA$ and $p\in\bN_{\geq 2}$
	that $\Pi_p(\cA)$ is a semiprime uniformly complete $f$-subalgebra of $\cA$, where the formulas
	(\ref{potenz}), (\ref{g1}) -- (\ref{f5a}) hold.
\par\medskip
We study now the finite elements in $\Pi_p(\cA)$.
\begin{theo}\label{T11}
Let $\cA$ be a uniformly complete $f$-algebra and let $p\in \bN_{\geq 2}$. Then   
\[g \in\Phi_1(\cA) \textit{ with a majorant } u \; \Lra\; g^p \in\Phi_1\big(\Pi_p(\cA)\big) \textit{ with the  majorant } u^p.\]
If, in addition, $\cA$ is semiprime, then
\[g \in\Phi_1(\cA) \textit{ with a majorant } u \; \Llra\; g^p \in\Phi_1\big(\Pi_p(\cA)\big) \textit{ with the  majorant } u^p.\]
\end{theo}
\prf
$\Ra$: Without loss of generality we assume $0<g\in\Phi_1(\cA)$. Otherwise consider $|g|$ and apply (${\rm d}_1$). 
If $u\in \cA_+$ is a majorant of $g$, then for each $f\in \cA_+$ there is a constant  
$c_{\hspace*{-.5mm}f}\geq 0$ with $f\Inf ng \leq c_{\hspace*{-.5mm}f} u$ for all $n\in\bN$. Then by means of formula (\ref{f5}) 
for the $p$-th power we get 
\be\label{f6} 
 f^p \Inf_p n^p g^p =f^p \Inf_p (ng)^p= (f\Inf ng)^p \leq (c_{\hspace*{-.5mm}f} u)^p = {c_{\hspace*{-.5mm}f}}^p u^p, \ee
where the last inequality follows\footnote{The twofold application of ($\ell_1$) on $0\leq a\leq b$ yields $a^2\leq b^2$. Indeed, by multiplying the inequality $0\leq a\leq b$ with $a$, resp. $b$, we obtain $a^2\leq ab$, resp. $ab\leq b^2$.} due to the condition ($\ell_1$).
\par
Let $f_1\cdots f_p \in\Pi_p^+(\cA)$ be an arbitrary element. By (\ref{g1}) we have $\Pi_p^+(\cA)=\Sigma_p(\cA)$, therefore there exists an $h\in \cA_+$ with $f_1\cdots f_p=h^p$.
Using (\ref{f6}) we get
\[h^p \Inf_p n^p g^p \leq {c_{\hspace{-0.2mm}h}}^p u^p \quad\mbox{for all } n\in\bN
\] 
and 
\[
(f_1\cdots f_p)\Inf_p n^p g^p \hs=\hs h^p \Inf_p n^p g^p \hs\leq\hs {c_{\hspace{-0.2mm}h}}^p u^p
          \quad \mbox{for all } n\in\bN.
\]
This shows that $g^p\in\Phi_1(\Pi_p(\cA))$ with the majorant $u^p$.
\par\smallskip
$\La$: Let $g^p$ be a positive finite element in $\Pi_p(\cA)$. 
There exist elements $u_1,\dots,u_p\in \cA_+$, such that for arbitrary $a_1,\dots,a_p\in \cA_+$ the inequality 
\be\label{f7}   
(a_1\cdots a_p)\Inf_p n g^p \hs\leq\hs c_{a_1\cdots a_p} (u_1\cdots u_p)  
\ee
holds for all $n\in\bN$ and some number $0<c_{a_1\cdots a_p}$. Since in $\cA$ there exists the element 
$u=(u_1\cdots u_p)^{\frac{1}{p}}$  the inequality (\ref{f7}) can be rewritten as 
\be\label{f8}
      (a_1\cdots a_p)\Inf_p n g^p \hs\leq\hs c_{a_1\cdots a_p} u^p. 
\ee
Now let be $a\in \cA_+$. 
By taking into consideration the relation (\ref{f5}) and the inequality (\ref{f8}) we get then  
\[ 
(a\Inf \sqrt[p]{n} \hs g)^p \hs=\hs a^p \Inf_p (\sqrt[p]{n} \hs g)^p \hs=\hs a^p \Inf_p n g^p\leq\hs c_{a^p}\,  u^p \quad \mbox{for all } n\in\bN.
\]
Due to the semiprimitivity the root is monotone and there holds
\[\big((a\Inf \sqrt[p]{n} \hs g)^p\big)^\frac{1}{p} \hs\leq\hs \sqrt[p]{c_{a^p}} \, u. \]
Therefore for all $n\in\bN$ there follows the inequality\footnote{For each $m\in \bN$ there exists $n\in \bN$ such that $m<\sqrt[p]{n}$, so that $a\Inf mg \hs\leq\hs a\Inf \sqrt[p]{n} g \hs\leq\hs \tilde{c}_a \,u$
for all $m\in\bN$.}
$ a\Inf \sqrt[p]{n} \hs g \hs\leq\hs \tilde{c}_{a} \,u \, \mbox{ with }\, \tilde{c}_{a}=\sqrt[p]{c_{a^p}}.  $
 This shows $g\in\Phi_1(\cA)$ with $u$ as one of its majorants.
\qed
\begin{cor}\label{C5}
Let $\cA$ be a uniformly complete $f$-algebra and let $p\in\bN_{\geq 2}$. Then
\begin{enumerate}
\item $ \; g_1,\dots,g_p$ are finite  in $\cA \quad\Lra\quad g_1\cdots g_p$ is finite in $\Pi_p(\cA)$.
\xdef\letzterwert{\the\value{enumi}}
\end{enumerate}
If, in addition, $\cA$ is semiprime, then
\begin{enumerate}
\setcounter{enumi}{\letzterwert}
\item $ \; g_1\cdots g_p$ is finite in $\Pi_p(\cA) \quad \Lra \quad (g_1\cdots g_p)^\frac{1}{p}$ 
             is finite  in $\cA$,
\item $ \; g_1,\dots,g_p$ are finite in $\cA \quad \Lra \quad (g_1\cdots g_p)^\frac{1}{p}$ 
             is finite in $\cA$,
\item $\;\Phi_1(\Pi_p(\cA))\subseteq \Phi_1(\cA)$,
\item $\Phi_1(\Pi_p(\cA)) = \Phi_1(\cA)\cap \Pi_p(\cA)$, \, provided $\cA$ has the weak factorization property.
\end{enumerate}
\end{cor}
\prf 
1. Let $g_1,\dots,g_p$ be positive finite elements in $\cA$. Then the element  
$g=g_1\Sup\dots\Sup g_p$ is also finite  in $\cA$ and, by the previous theorem the element $g^p$ is finite in $\Pi_p(A)$.
Since $0\leq g_i\leq g$ for all $i=1,\ldots, p$ by condition ($\ell_1$)  we have 
\[g_1\cdots g_p \hs\leq\hs g_1\cdots g_{p-1}g \hs\leq\hs g_1\cdots g_{p-2}g^2
		\hs\leq\hs \dots \hs\leq\hs g_1 g^{p-1} \hs\leq\hs g^p.  \]
The element $g_1\cdots g_p$ is finite in $\Pi_p(\cA)$ since $\Phi_1\big(\Pi_p(\cA)\big)$ is an order ideal in 
$\Pi_p(\cA)$.
\par
Let $g_1,\dots,g_p$ be arbitrary finite elements in $\cA$.
By the first part of the proof the element $|g_1|\cdots|g_p|$ is finite  in $\Pi_p(\cA)$.
Due to (\ref{f5a}) we have $|g_1|\cdots|g_p|=|g_1\cdots g_p|_p$ and so $g_1\cdots g_p$ is a finite element in $\Pi_p(\cA)$.\\
2. Without loss of generality we may assume that $g_1,\dots,g_p$ are positive elements in $\cA$, otherwise consider $|g_1|,\dots,|g_p|$ and apply (\ref{f5a}). 
According to Remark 3(1) there exists an element $\tilde{g}=(g_1\cdots g_p)^\frac{1}{p}$ in $\cA$.
The equality $g_1\cdots g_p = ((g_1\cdots g_p)^{\frac{1}{p}})^p = \tilde{g}^p$ 
shows that the finiteness of $g_1\cdots g_p$ in $\Pi_p(\cA)$ implies that $\tilde{g}^p$ is finite. 
By the  theorem one has $\tilde{g}\in\Phi_1(\cA)$. \\
3.    Follows directly from 1. and 2.\\
4.    Let $g_1\cdots g_p \in \Phi_1(\Pi_p(\cA))$. Then by part 2. we get
$(g_1\cdots g_p)^{\frac{1}{p}} \in\Phi_1(\cA)$, which according to Theorem \ref{T11} yields $g_1\cdots g_p \in \Phi_1(\cA)$.\\
5.    The relation "$\subseteq$" follows from 4. For the converse relation "$\supseteq$" let $\varphi\in\Phi_1(\cA)\cap\Pi_p(\cA)$. Then the element $\varphi$ can be written as a $p$-fold product $\varphi=\varphi_1\dots \varphi_p$ and therefore possesses the root $\varphi^{\frac{1}{p}} \in\cA$. By Theorem~\ref{T8} we have $\varphi\in\Phi_1(\cA)$ and by means of 1. then $\varphi\in\Phi_1(\Pi_p(\cA))$.\qed
\par\medskip
In the next corollary we obtain some information on the totally finite and the self-majorizing elements in an $f$-algebra.
For its proof we need the following 
\begin{Lem}\label{L2}
Let  $\cA$ be a uniformly complete $f$-algebra and let $p\in \bN_{\geq 2}$. Then for all $g\in\cA$ the following implication holds:
\be\label{S(A)thenS(Pi(A))} g\in S(\cA) \quad\Lra\quad g^p\in S\big(\Pi_p(\cA)\big).\ee
If, in addition, $\cA$ is semiprime, then
\be\label{S(A)iffS(Pi(A))} g\in S(\cA) \quad\Llra\quad g^p\in S\big(\Pi_p(\cA)\big).\ee
\end{Lem}
\prf 
Let $g$ be a self-majorizing element in $\cA$, i.e. $|g|$ is a majorant of $g$ in $\cA$.
By Theorem~\ref{T11} this implies that $|g|^p$ is a majorant of $g^p$ in $\Pi_p(\cA)$.
Formula (\ref{f5a}) yields 
the equality $|g|^p = |g^p|_p$, so $|g^p|_p\,$ is a majorant of $g^p$ in $\Pi_p(\cA)$. 
Therefore $g^p \in S\big(\Pi_p(\cA)\big)$.
\par
Conversely, let $g^p\in S\big(\Pi_p(\cA)\big)$, i.e. $|g^p|_p$ is a majorant of $g^p$ in $\Pi_p(\cA)$.
The equa\-li\-ty $|g^p|_p = |g|^p$ and Theorem~\ref{T11} imply that $|g|$ is a majorant of $g$ in 
$\cA$ and therefore, $g\in S(\cA)$.
\qed 
\begin{cor}\label{C6}
Let $\cA$ be a uniformly complete $f$-algebra and let $p\in\bN_{\geq2}$.
Then for all $g\in \cA$ the following implications hold:
\begin{enumerate}
\item $g\in\Phi_2(\cA) \;\Lra\; g^p \in\Phi_2\big(\Pi_p(\cA)\big)$,
\item $g\in\Phi_3(\cA) \;\Lra\; g^p \in\Phi_3\big(\Pi_p(\cA)\big)$.
\end{enumerate}
If, in addition, $\cA$ is semiprime, then the converse implications are also true.
\end{cor}
\prf 1.  
$\Rightarrow$: Let $g\in\Phi_2(\cA)$ have a finite majorant $u\in\cA$.
By the first part of Theorem~\ref{T11} 
we obtain $g^p\in \Phi_1\big(\Pi_p(\cA)\big)$ with majorant $u^p$, and   
the same theorem guarantees the finiteness of the majorant $u^p$ in $\Pi_p(\cA)$, i.e. $g^p\in \Phi_2\big(\Pi_p(\cA)\big)$.
\par
$\Leftarrow$: Let $g^p\in\Phi_2\big(\Pi_p(\cA)\big)$ with a finite majorant $u_1\cdots u_p$.
By Remark 3(1) we can write this majorant as a $p$-fold product
$u_1\cdots u_p= \big((u_1\cdots u_p)^\frac{1}{p}\big)^p = u^p$ of the element $u= (u_1\cdots u_p)^\frac{1}{p}$.
Then the semiprimitivity of $\cA$ and Theorem~\ref{T11} yield $g\in \Phi_1(\cA)$ with the majorant $u$ and also 
the finiteness of the majorant $u$ in $\cA$. Therefore $g\in \Phi_2(\cA)$. 
\par\smallskip
2. The set $\Phi_3(\cA)$ coincides with the order ideal generated by the set $S(\cA)$ (see \cite[Corollary 2]{TeiWeb11}), i.e.
\[\Phi_3(\cA)= \{a\in \cA\colon \exists s_1,\ldots,s_n\in S(\cA)\;\mbox{and}\;\la_1,\ldots,\la_n\in \bR_{\geq 0} \mbox{ with }  \bl a\br\leq \sum_{i=1}^n\la_i\bl s_i\br \}.\]
Since  $\sum\limits_{i=1}^n\la_i\bl s_i\br$ is a positive self-majorizing element (see \cite[Proposition 1]{TeiWeb11}), the order ideal 
$\Phi_3(\cA)$ can be written as $\Phi_3(\cA) = \left\{a\in\cA\colon \exists s\in S_+(\cA): |a|\leq s\right\}$.
\par
$\Rightarrow$:
Let $g\in\Phi_3(\cA)$.  
There is an $s\in S_+(\cA)$ such that $|g|\leq s$.
Since $s$ is a majorant of $s$ in $\cA$, the element $s$ is also a majorant of $|g|$.
By the first part of Theorem~\ref{T11} we obtain that the element $|g|^p$ is finite in $\Pi_p(\cA)$ with a majorant $s^p$.
Due to (\ref{S(A)thenS(Pi(A))}) the element $s^p$ is self-majorizing in $\Pi_p(\cA)$.
The formula (\ref{f5a}) yields 
that the element $|g^p|_p$ belongs to the order ideal generated by $S_+\big(\Pi_p(\cA)\big)$,
i.e. $g^p \in \Phi_3\big(\Pi_p(\cA)\big)$.
\par
$\Leftarrow$:
Conversely, let $g^p \in\Phi_3\big(\Pi_p(\cA)\big)$.
Since $\Phi_3\big(\Pi_p(\cA)\big)$ is the order ideal generated by $S_+\big(\Pi_p(\cA)\big)$ in $\Pi_p(\cA)$,
there is an element $s\in S_+\big(\Pi_p(\cA)\big)$ such that $|g^p|_p\leq s$.
Using Remark 3(1) we can write the majorant $s$ as $s = s_1\cdots s_p = \tilde{s}^{\;p}$,
where $\tilde{s}:= (s_1\cdots s_p)^\frac{1}{p}$.
Notice that $\tilde{s}^{\;p}$ has itself as a majorant in $\Pi_p(\cA)$.
Due to (\ref{S(A)iffS(Pi(A))}) and the second part of Theorem~\ref{T11} the element $\tilde{s}$ is self-majorizing
in $\cA$ and is a majorant of $g$ in $\cA$. Therefore we obtain $g\in\Phi_3(\cA)$.
\qed
\par\medskip
By summing up the results obtained in Theorem~\ref{T11}, 
Corollaries~\ref{C5} and \ref{C6} we may write
\begin{cor}\label{C7}
Let  $\cA$ be a semiprime uniformly complete $f$-algebra and $p\in \bN_{\geq 2}$.
Then for $i=1,2,3$ there holds
\[      \big(\Phi_i(\cA)\big)^p = \Phi_i\big(\Pi_p(\cA)\big),    \]
where
$\big(\Phi_i(\cA)\big)^p \hs=\hs \left\{ g_1\cdots g_p\in\Pi_p(\cA) 
  \colon  g_1,\dots,g_p \in\Phi_i(\cA)\right\}$.
\end{cor}
\prf
Let $i=1$. Indeed, the relation "$\subseteq$" follows from Corollary \ref{C5}.\hspace*{.5mm}1. 
The relation "$\supseteq$" is obtained as follows: 
Let $g\in\Phi_1(\Pi_p(\cA))$, i.e. $g=g_1\cdots g_p$ with $g_j\in \cA$, \mbox{$j=1,\ldots,p$}. 
Then by Corollary \ref{C5}.\hspace*{.5mm}2. the element $g^{\frac{1}{p}}$ is finite in $\cA$.
From $g = (g^{\frac{1}{p}})^p$ it is clear that $g$ is a product consisting of $p$ finite elements of $\cA$, i.e. $g\in\big(\Phi_1(\cA)\big)^p$.
\par
The cases $i=2,3$ are proved similarly using Corollary~\ref{C6}.
\qed
\par\medskip
The proof of the second inclusion of Corollary \ref{C7} (for $i=1$) shows 
that each finite element of $\Pi_p(\cA)$ has a representation as the $p$-th power of a single finite element of $\cA$.
In general, $g = g_1 \cdots g_p \in \Phi_1\big(\Pi_p(\cA)\big)$ does not imply 
$g_1,\dots,g_p \in \Phi_1(\cA)$, what is demonstrated by the next example. 
\par\smallskip
\hspace*{1cm}\\
\bsp{\bf 4}\;\label{ex4} 
Let $\cA=C([0,\infty))$ be the vector lattice of  all continuous functions on the interval $[0, \infty)$ and equip $\cA$ with the pointwise order and the algebraic operations. Then $\cA$ is an Archimedean unitary semiprime uniformly complete  $f$-algebra. For $p=3$ consider  
$\Pi_3(\cA) = \left\{ f_1 f_2 f_3 \colon f_1, f_2, f_3 \in \cA\right\}$.
Since the function \ou$_{[0,\infty)}$ is the multiplicative unit in $\cA$, all functions of  $\cA$ belong to $\Pi_3(\cA)$. 
This means $\cA$ and $\Pi_3(\cA)$ coincide. 
\par
The finite elements in $\cA$ are exactly the functions with compact support. 
Consider the following three functions of $\cA$:
\[ f_1(t) = t, \qquad f_2(t) = \mbox{\ou}_{[0,\infty)}
		\quad \mbox{ and } \quad f_3(t)=
                 \left\{\begin{array}{ll} 
			    \sin t & \mbox{for }t\in [0,\pi], \\
			                    0 & \mbox{for }t\in (\pi,\infty).
       				          \end{array}\right. 
\]
The only finite element among them is $f_3$. 
The product $f_1 f_2 f_3$, i.e. the function
\[\varphi(t)=\left\{\begin{array}{ll} t \sin t & \mbox{for }t\in [0,\pi], \\
  0 & \mbox{for }t\in (\pi, \infty),\end{array}\right. \]
is a finite element in $\cA=\Pi_3(\cA)$, however not all of its factors are finite elements. 
\par
In view of Corollary \ref{C7} we know that there exists a finite function  $\tilde{\varphi}$ in  
$\cA$ such that $\tilde{\varphi}^3 = f_1 f_2 f_3$. In our case this is the function
\par
\vspace*{0.5cm}
\begin{minipage}[]{.4\textwidth} 
\[
\tilde{\varphi}(t) = \left\{\begin{array}{ll} (t\sin t)^{\frac{1}{3}} & \mbox{for } t\in [0,\pi], \\
        				0 & \mbox{for } t\in (\pi, \infty).\end{array}\right.
\]\\ \\ \\
\end{minipage}
\hspace*{.1\textwidth}
\begin{minipage}[]{.4\textwidth}
\includegraphics[]{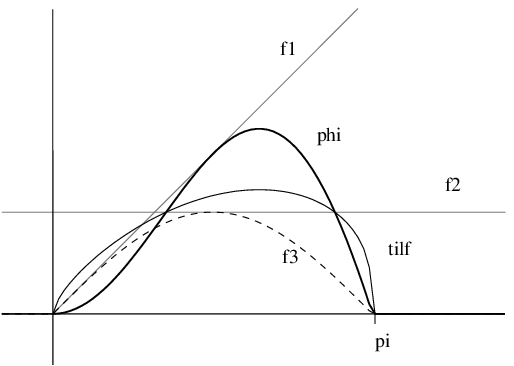}
\end{minipage}
\par\bigskip
\bigskip
\textbf{Acknowledgments:} The authors thank the referee for valuable remarks and suggestions, which enabled us to present our results more transparently and precisely.
\bibliographystyle{plain}

\begin{thebibliography}{00}
%
\bibitem{AbrAli02}
Y.A. Abramovich, C.D. Aliprantis:
{An Invitation to Operator Theory.} 
{Graduate Studies in Mathematics, vol.50.}
\newblock Amer.Math.Soc., Providence, Rhode Island (2002) 
%
\bibitem{AliBur85}
 C.D.~Aliprantis, O.~Burkinshaw. 
\newblock {Positive Operators.} 
\newblock Academic Press{,} Inc., London (1985)
%
\bibitem{BerHui90}
S.J.~Bernau, C.B.~Huijsmans.
\newblock Almost $f$-algebras and $d$-algebras.
\newblock {Math. Proc.Cambridge Philos. Soc.}, {\bf{107}}(2), 287--308 (1990)
%
\bibitem{BeuHui84}
F.~Beukers and C.B.~Huijsmans.
\newblock Calculus in $f$-algebras.
\newblock {J. Austral. Math. Soc.}(Series A), {\bf{37}}(1), 110--116 (1984)
%
\bibitem{BeuHuideP83}
F.~Beukers, C.B.~Huijsmans and B.~de~Pagter.
\newblock Unital embedding and complexification of $f$-algebras.
\newblock {Math. Z.}, {\bf{183}}(1), 131--144 (1983)
%
\bibitem{BKW77}
 A.~Bigard, K.~Keimel, S.~Wolfenstein. 
\newblock {Groupes ez Anneaux R$\acute{e}$ticul$\acute{e}$s.} 
           Lecture Notes in Mathematics 608. 
\newblock Springer-Verlag{,} Berlin{,} Heidelberg{,} New York (1977)
%
\bibitem{Bou00}
K.~Boulabiar.
\newblock Products in almost $f$-algebras.
\newblock {Comment. Math. Univ. Carolinae.}, {\bf{41}}, 747--759 (2000)
%
\bibitem{Bou03}
K.~Boulabiar.
\newblock On products in lattice-ordered algebras.
\newblock {J. Austral. Math. Soc.}, {\bf{75}}, 23--40 (2003)
%
\bibitem{BBT03}
K.~Boulabiar, G.~Buskes, A.~Triki.
\newblock Some Recent Trends and Advances in certain Lattice Ordered Algebras.
\newblock {Contemp. Math.} 328, 99--133, Amer. Math. Soc., Providence (2003)
%
\bibitem{ChenWeb03-1}
Z.L.~Chen and M.R.~Weber.
\newblock  On finite elements in vector lattices and {B}anach lattices.
\newblock {Math. Nachrichten}, {\bf{279}}(5-6), 495--501 (2006)
%
\bibitem{ChenWeb03-2}
Z.L.~Chen and M.R.~Weber.
\newblock On finite elements in sublattices of {B}anach lattices.
\newblock {Math. Nachrichten}, {\bf{280}}(5-6), 485--494 (2007) 
%
\bibitem{ChenWeb03-3}
Z.L.~Chen and M.R.~Weber.
\newblock On finite elements in lattices of regular operators.
\newblock {Positivity}, {\bf{11}}, 563--574 (2007)
%
\bibitem{deP81}
B.~de Pagter. 
\newblock {$f$-Algebras and Orthomorphisms.} 
\newblock Thesis, Leiden (1981)
%
\bibitem{FelPor76}
W.A.~Feldman and J.F.~Porter.
\newblock Order units and base norms generalized for convex spaces.
\newblock {\em Proc. London Math. Soc.}, Vol. XXXIII, Sept. 1976.
%
\bibitem{HHW} 
N.~Hahn, S.~Hahn and M.R.~Weber.
\newblock On finite elements in vector lattices of operators.   
\newblock {Positivity}, Online, {\bf{13}}(1), 145--163 (2009) 
%
\bibitem{Hui91}
C.B.~Huismans.
\newblock Lattice-ordered algebras and $f$-algebras: a survey.
\newblock {Proceedings of the Conference: {\rm Positive Operators, Riesz Spaces, and Economics}}, 151--169, Eds.: C.D.~Aliprantis, K.C.~Border, W.A.J.~Luxemburg. Caltech, Pasadena (1991)
%
\bibitem{Kud62}
V.~Kudl\'{a}\v{c}ek.
\newblock {O n\v{e}kter\'{y}ch typech $\ell$-okruhu. (Czech)}. 
\newblock  Sbornik Vysok\'{e}ho U\v{c}eni Tech. v Brn\v{e}{,} {\bf{1-2}}{,} 179-181 (1962)
%
\bibitem{LuxMoo67}
W.A.J.~Luxemburg and L.C.~Moore, Jr.
\newblock Archimedean quotient Riesz spaces.
\newblock {\em Duke Math. J.}, 34:725--739, 1967
%
\bibitem{LuxZaa71}
W.A.J.~Luxemburg, A.C.~Zaanen.
\newblock {Riesz Spaces I}.
\newblock  North-Holland Publ. Comp.{,} Amsterdam (1971)
%
\bibitem{MakWeb74}
B.M.~Makarow and M.~Weber.
\newblock {\"U}ber die {R}ealisierung von {V}ektorverb{\"a}nden {I}.
  ({R}ussian).
\newblock {Math. Nachrichten}, {\bf{60}}, 281--296 (1974)
%
\bibitem{MakWeb77a}
B.M.~Makarow and M.~Weber.
\newblock Einige {U}ntersuchungen des {R}aumes der maximalen {I}deale eines
  {V}ektorverbandes mit {H}ilfe finiter {E}lemente {I}.
\newblock {Math. Nachrichten}, {\bf{79}}, 115--130 (1977)
%
\bibitem{MakWeb77b}
B.M.~Makarow and M.~Weber.
\newblock Einige {U}ntersuchungen des {R}aumes der maximalen {I}deale eines
  {V}ektorverbandes mit {H}ilfe finiter {E}lemente {II}.
\newblock {Math. Nachrichten}, {\bf{80}}, 115--125 (1977)
%
\bibitem{MakWeb78}
B.M.~Makarow and M.~Weber.
\newblock {\"U}ber die {R}ealisierung von {V}ektorverb{\"a}nden {III}.
\newblock {Math. Nachrichten}, {\bf{86}}, 7--14 (1978)
%
\bibitem{Mey91}
P.~Meyer-Nieberg.
\newblock {Banach Lattices}.
\newblock Springer-Verlag, Berlin{,} Heidelberg{,} New York (1991)
%
\bibitem{TeiWeb11} 
K.~Teichert and M.R.~Weber.
\newblock On self-majorizing elements in Archimedean vector lattices. 
\newblock {Preprint}, MATH-AN-01-2011, TU Dresden (2011)
%
\bibitem{Web73}
M.~Weber.
\newblock {O}n the representation of vector lattices with a countable fundamental sequence of intervals. (Russian),
\newblock {Vestnik Leningrad Univ.}, {\bf{3}}(7), 152--154 (1973)
%
\bibitem{Web95b}
M.R.~Weber.
\newblock On finite and totally finite elements in vector lattices.
\newblock {Analysis Mathematica}, {\bf{21}}, 237--244 (1995)
%
\bibitem{Web06}
M.R.~Weber.
\newblock Finite elements in vector lattices.
\newblock {Proceedings of  the Conference Positivity IV}, 155--172, Eds.: J.~Voigt, M.R.~Weber. Dresden (2006)
%
\bibitem{Zaa83}
A.C.~Zaanen,
\newblock {Riesz Spaces II}. 
\newblock  North-Holland Publ. Comp.{,} Amsterdam (1983)
%
\bibitem{Zaa97}
A.C.~Zaanen.
\newblock {Introduction to Operator Theory in Riesz Spaces}.
\newblock Springer-Verlag, Berlin{,} Heidelberg{,} New York (1997)
%
\end{thebibliography}

\end{document}